\documentclass[11pt]{amsart}
\usepackage{amssymb}
\usepackage{amsmath}
\usepackage{fancyhdr}
\usepackage[british]{babel}
\usepackage{geometry}
\usepackage{enumitem}
\usepackage{algpseudocode}
\usepackage{dsfont}
\usepackage{centernot}
\usepackage{xstring}
\usepackage{colortbl}
\usepackage[section]{algorithm}
\usepackage{graphicx}
\usepackage[font={footnotesize}]{caption}
\usepackage[usenames,dvipsnames,table]{xcolor}
\usepackage{pstricks,tikz}
\usepackage[h]{esvect}
\usepackage[unicode,
		bookmarksopen=true,
		bookmarksopenlevel=1,
		colorlinks=true,
		linkcolor=darkblue,
        linktoc=page,
		citecolor=darkblue,
]{hyperref}

\title[]{Triple systems with no three triples spanning at most five points}
\author[S.~Glock]{Stefan Glock}

\address{School of Mathematics, University of Birmingham,
Edgbaston, Birmingham, B15 2TT, United Kingdom}
\email{s.glock@bham.ac.uk}

\thanks{The research leading to these results was partially supported by the European Research Council
under the European Union's Seventh Framework Programme (FP/2007--2013) / ERC Grant
Agreement no. 306349.}

\newtheorem{theorem}{Theorem}

\newtheorem{lemma}[theorem]{Lemma}

\newtheorem{conj}[theorem]{Conjecture}

\theoremstyle{definition}

\newtheoremstyle{claimstyle}{5pt}{5pt}{\em}{5pt}{\em}{:}{5pt}{}
\theoremstyle{claimstyle}

\newtheoremstyle{stepstyle}{10pt}{5pt}{\em}{0pt}{\em}{:}{5pt}{}
\theoremstyle{stepstyle}

\definecolor{darkblue}{rgb}{0,0,0.5}

\def\noproof{{\unskip\nobreak\hfill\penalty50\hskip2em\hbox{}\nobreak\hfill%
       $\square$\parfillskip=0pt\finalhyphendemerits=0\par}\goodbreak}
\def\endproof{\noproof\bigskip}

\newdimen\margin
\def\textno#1&#2\par{
   \margin=\hsize
   \advance\margin by -4\parindent
          \setbox1=\hbox{\sl#1}
   \ifdim\wd1 < \margin
      $$\box1\eqno#2$$
   \else
      \bigbreak
      \hbox to \hsize{\indent$\vcenter{\advance\hsize by -3\parindent
      \it\noindent#1}\hfil#2$}
      \bigbreak
   \fi}

\def\proof{\removelastskip\penalty55\medskip\noindent\setcounter{claim}{0}\setcounter{step}{0}{\bf Proof. }} % in each main proof, claim counter set back

\begin{document}

\def\COMMENT#1{}
\def\TASK#1{}

\def\eps{{\varepsilon}}
\newcommand{\ex}{\mathbb{E}}
\newcommand{\pr}{\mathbb{P}}
\newcommand{\cB}{\mathcal{B}}
\newcommand{\cA}{\mathcal{A}}
\newcommand{\cE}{\mathcal{E}}
\newcommand{\cS}{\mathcal{S}}
\newcommand{\cF}{\mathcal{F}}
\newcommand{\cG}{\mathcal{G}}
\newcommand{\bL}{\mathbb{L}}
\newcommand{\bF}{\mathbb{F}}
\newcommand{\bZ}{\mathbb{Z}}
\newcommand{\cH}{\mathcal{H}}
\newcommand{\cC}{\mathcal{C}}
\newcommand{\cM}{\mathcal{M}}
\newcommand{\bN}{\mathbb{N}}
\newcommand{\bR}{\mathbb{R}}
\def\O{\mathcal{O}}
\newcommand{\cP}{\mathcal{P}}
\newcommand{\cQ}{\mathcal{Q}}
\newcommand{\cR}{\mathcal{R}}
\newcommand{\cJ}{\mathcal{J}}
\newcommand{\cL}{\mathcal{L}}
\newcommand{\cK}{\mathcal{K}}
\newcommand{\cD}{\mathcal{D}}
\newcommand{\cI}{\mathcal{I}}
\newcommand{\cV}{\mathcal{V}}
\newcommand{\cT}{\mathcal{T}}
\newcommand{\cU}{\mathcal{U}}
\newcommand{\cX}{\mathcal{X}}
\newcommand{\cZ}{\mathcal{Z}}
\newcommand{\1}{{\bf 1}_{n\not\equiv \delta}}
\newcommand{\eul}{{\rm e}}
\newcommand{\Erd}{Erd\H{o}s}
\newcommand{\cupdot}{\mathbin{\mathaccent\cdot\cup}}
\newcommand{\whp}{whp }

\newcommand{\doublesquig}{%
  \mathrel{%
    \vcenter{\offinterlineskip
      \ialign{##\cr$\rightsquigarrow$\cr\noalign{\kern-1.5pt}$\rightsquigarrow$\cr}%
    }%
  }%
}

\newcommand{\defn}{\emph}

\newcommand\restrict[1]{\raisebox{-.5ex}{$|$}_{#1}}

\newcommand{\prob}[1]{\mathrm{\mathbb{P}}\left(#1\right)}
\newcommand{\cprob}[2]{\mathrm{\mathbb{P}}_{#1}\left(#2\right)}
\newcommand{\expn}[1]{\mathrm{\mathbb{E}}\left(#1\right)}
\newcommand{\cexpn}[2]{\mathrm{\mathbb{E}}_{#1}\left(#2\right)}
\def\gnp{G_{n,p}}
\def\G{\mathcal{G}}
\def\lflr{\left\lfloor}
\def\rflr{\right\rfloor}
\def\lcl{\left\lceil}
\def\rcl{\right\rceil}

\newcommand{\qbinom}[2]{\binom{#1}{#2}_{\!q}}
\newcommand{\binomdim}[2]{\binom{#1}{#2}_{\!\dim}}

\newcommand{\grass}{\mathrm{Gr}}

\newcommand{\brackets}[1]{\left(#1\right)}
\def\sm{\setminus}
\newcommand{\Set}[1]{\{#1\}}
\newcommand{\set}[2]{\{#1\,:\;#2\}}
\newcommand{\krq}[2]{K^{(#1)}_{#2}}
\newcommand{\ind}[1]{$(\ast)_{#1}$}
\newcommand{\indcov}[1]{$(\#)_{#1}$}
\def\In{\subseteq}

\begin{abstract}  \noindent
We show that the maximum number of triples on $n$~points, if no three triples span at most five points, is $(1\pm o(1))n^2/5$. More generally, let $f^{(r)}(n;k,s)$ be the maximum number of edges of an $r$-uniform hypergraph on $n$~vertices not containing a subgraph with $k$~vertices and $s$~edges. In 1973, Brown, Erd\H{o}s and S\'os conjectured that the limit $\lim_{n\to \infty}n^{-2}f^{(3)}(n;k,k-2)$ exists for all~$k$. They proved this for $k=4$, where the limit is $1/6$ and the extremal examples are Steiner triple systems. We prove the conjecture for $k=5$ and show that the limit is~$1/5$. The upper bound is established via a simple optimisation problem. For the lower bound, we use approximate $H$-decompositions of~$K_n$ for a suitably defined graph~$H$.
\end{abstract}

\maketitle

\section{Introduction}

For a family $\cF$ of $r$-graphs (i.e.~$r$-uniform hypergraphs), let $ex(n;\cF)$ denote the maximum number of edges in an $\cF$-free $r$-graph on $n$~vertices, which is called the \defn{Tur\'an number of $\cF$}. %One of the most famous open problems in Extremal Combinatorics is his conjecture that $ex(n;K_4^{(3)})=(\frac{5}{9}\pm o(1))\binom{n}{3}$.
Here, we consider the family $\cF^{(r)}(k,s)$ of all $r$-graphs on $k$ vertices with $s$~edges.
In 1973, Brown, Erd\H{o}s and S\'os introduced the function $f^{(r)}(n;k,s):=ex(n;\cF^{(r)}(k,s))$. A lot of research has been invested to understand this function asymptotically (e.g.~\cite{AS:06,BES:73a,BES:73b,erdos:64,EFR:86,GS:17,RS:78,SS:05}).
Using the probabilistic method, Brown, Erd\H{o}s and S\'os~\cite{BES:73b} showed that $f^{(r)}(n;k,s)=\Omega\left(n^{(rs-k)/(s-1)}\right)$ for all $k>r$ and $s\ge 2$. They deduced that
$f^{(3)}(n;k,k-2)=\Theta(n^2)$ for every fixed $k\ge 4$, and posed the following conjecture.
\begin{conj}[Brown, Erd\H{o}s, S\'os~\cite{BES:73b}] \label{conj:BES}
The limit $\lim_{n\to \infty}n^{-2}f^{(3)}(n;k,k-2)$ exists for all~$k\ge 4$.
\end{conj}
They confirmed this for $k=4$, where the limit is~$1/6$. For $k=5$, they gave a lower bound of~$1/6$ and an upper bound of~$2/9$. Here, we prove the conjecture for $k=5$ and show that the limit is~$1/5$.

\begin{theorem} \label{thm:main}
$\lim_{n\to \infty}n^{-2}f^{(3)}(n;5,3)=\frac{1}{5}$.
\end{theorem}

We believe that our methods can lead to further progress concerning Conjecture~\ref{conj:BES} and related questions.

%Previously, the best known lower bound was~$1/6$, and the best known upper bound~$2/9$.

%Note that for $s=\binom{k}{r}$, we have $\cF^{(r)}(k,s)=\Set{K^{(r)}_k}$. Tur\'an conjectured that $\lim_{n\to \infty}n^{-3}f^{(3)}(n;4,4)=5/9$.

%$f^{(r)}(n;s(r-i)+i,s)=\Theta(n^i)$ for all $2\le i<r$ and $s\ge 2$. In fact, they showed that for $k>r$ and $s\ge 2$, $f^{(r)}(n;k,s)=\Omega(n^\frac{rs-k}{s-1})$.

\subsection{Results for \texorpdfstring{$f^{(3)}(n;k,k-2)$}{f(n;k,k-2)}}

%Define $$c_k^-:=\liminf_{n\to \infty}n^{-2}f^{(3)}(n;k,k-2) \quad \mbox{and}\quad c_k^+:=\limsup_{n\to \infty}n^{-2}f^{(3)}(n;k,k-2).$$
For a $3$-graph $G$ and a pair $x,y$ of distinct vertices, we let $d(xy)$ denote the \defn{codegree} of $xy$, that is, the number of edges containing $x$ and~$y$. We call $G$ \defn{linear} if the maximum codegree is at most~$1$. A \defn{Steiner triple system of order $n$} is a $3$-graph on $n$ vertices such that all codegrees are equal to~$1$. Due to an old theorem of Kirkman, such systems exist if and only if $n\equiv 1,3\mod{6}$.

Brown, Erd\H{o}s and S\'os~\cite{BES:73a} showed that $\lim_{n\to \infty}n^{-2}f^{(3)}(n;4,2)=1/6$. Clearly, any $3$-graph on $n$ vertices with more than $\binom{n}{2}/3$ edges contains a pair of vertices with codegree at least~$2$, and is thus not $\cF^{(3)}(4,2)$-free. On the other hand, any Steiner triple system is $\cF^{(3)}(4,2)$-free and has $\binom{n}{2}/3$ edges.

Moreover, since any Steiner triple system is also $\cF^{(3)}(5,3)$-free, this yields the mentioned lower bound $\liminf_{n\to \infty} n^{-2}f^{(3)}(n;5,3) \ge 1/6$.
Perhaps this led Erd\H{o}s~\cite{erdos:73,erdos:76} to his conjecture on the existence of locally sparse Steiner triple systems. More precisely, he conjectured that for any~$k$, there exists $n_0$ such that for all $n\ge n_0$, there exists a Steiner triple system of order~$n$ which is $\bigcup_{4\le j\le k}\cF^{(3)}(j,j-2)$-free (subject to the necessary condition $n\equiv 1,3\mod{6}$). Such Steiner triple systems are also referred to as having large `girth'.
This conjecture was proved asymptotically in~\cite{GKLO:18}, and independently in~\cite{BW:18}, by showing that for any fixed~$k$, as $n\to \infty$, there exists a $\bigcup_{4\le j\le k}\cF^{(3)}(j,j-2)$-free $3$-graph $G$ on $n$ vertices with $(1/6-o(1))n^2$ edges. %In both papers, this is shown by analysing a natural random process.
In particular, this implies that for every~$k\ge 4$,
\begin{align}
\liminf_{n\to \infty} n^{-2}f^{(3)}(n;k,k-2) \ge \frac{1}{6}, \label{best general lower}
\end{align}
which is to date the best lower bound for Conjecture~\ref{conj:BES}.
When only considering linear $3$-graphs, this would be best possible. Moreover, Steiner triple systems are maximal in the sense that adding any further triple creates a forbidden subgraph. %with $k-2$ edges on at most $k$~vertices for every $k\ge 4$.
In view of this, one may ask whether \eqref{best general lower} is sharp in general, or whether we can pack significantly more edges into an $\cF^{(3)}(k,k-2)$-free $3$-graph $G$ if we do not require $G$ to be linear.
%Prior to this paper, no such example was known.

This leads to the discussion of upper bounds. A trivial upper bound is given by $f^{(3)}(n;k,k-2)\le (k-3)\binom{n}{2}/3$. Indeed, any $3$-graph on $n$ vertices with more than $(k-3)\binom{n}{2}/3$ edges contains a pair of vertices with codegree at least~$k-2$, and is thus not $\cF^{(3)}(k,k-2)$-free.
As indicated in~\cite{BES:73a}, this can be improved significantly to
\begin{align*}
f^{(3)}(n;k,k-2)\le \frac{k-3}{3(k-2)}n(n-1), %\label{best general upper}
\end{align*}
by averaging over vertex degrees instead of codegrees, and using the fact that $f^{(2)}(n;k,k-1)=\lflr\frac{k-2}{k-1}n\rflr$~(see~\cite{erdos:64}). Indeed, if $G$ is a $3$-graph on $n$ vertices with $e(G)>\frac{k-3}{3(k-2)}n(n-1)$, then some vertex $x$ has degree larger than $\frac{k-3}{k-2}(n-1)\ge f^{(2)}(n-1;k-1,k-2)$. This yields an $\cF^{(2)}(k-1,k-2)$-graph in the link graph of~$x$, and thus an $\cF^{(3)}(k,k-2)$-graph in~$G$.

To sum up, the currently best known bounds for Conjecture~\ref{conj:BES} are
\begin{align}
\frac{1}{6} \le \liminf_{n\to \infty} n^{-2}f^{(3)}(n;k,k-2) \le \limsup_{n\to \infty} n^{-2}f^{(3)}(n;k,k-2) \le \frac{k-3}{3(k-2)}. \label{best general}
\end{align}

For $k=5$, we show that neither of these bounds gives the correct answer, and there is not much reason to believe that this changes for larger~$k$. Brown, Erd\H{o}s and S\'os~\cite{BES:73a} actually suggested that the correct answer for $k=6$ should be~$1/4$, matching the upper bound in~\eqref{best general}. However, our methods can be easily used to refute this (cf.~Section~\ref{sec:further}).

\subsection{Results for \texorpdfstring{$f^{(3)}(n;k,k-3)$}{f(n;k,k-3)}}
%Whilst for the $(k,k-2)$-problem, at least the asymptotics are known, much less is known for the $(k,k-3)$-problem, which we discuss briefly.
The above conjecture of Erd\H{o}s is best possible in the sense that every Steiner triple system of order~$n$ contains an $\cF^{(3)}(k,k-3)$-graph for every $4\le k\le n$. This is true in a very robust sense. For instance, Ruzsa and Szemer\'edi~\cite{RS:78} showed that $n^{2-o(1)}<f^{(3)}(n;6,3)=o(n^2)$, which solved a problem of Brown, Erd\H{o}s and S\'os~\cite{BES:73a,BES:73b} and has become known as the \defn{$(6,3)$-theorem}. The $(6,3)$-theorem is closely related with the development of the regularity lemma and the triangle removal lemma, and bounds for Roth’s theorem. Moreover, the problem can be translated into an induced matching problem in graphs (see also~\cite{FHS:17}).
Erd\H{o}s, Frankl, and R{\"o}dl~\cite{EFR:86} extended this result to any~$r$, showing that $n^{2-o(1)} < f^{(r)}(n;3(r-2)+3,3)=o(n^2)$. Alon and Shapira~\cite{AS:06} extended this result further by showing that $n^{j-o(1)} < f^{(r)}(n;3(r-j)+j+1,3)=o(n^j)$ for any $r>j\ge 2$, and also generalised a conjecture from~\cite{EFR:86} to the following.

\begin{conj}[cf.~\cite{AS:06}]
For any $r>j\ge 2$ and $s\ge 3$, we have $$n^{j-o(1)} < f^{(r)}(n;s(r-j)+j+1,s)=o(n^j).$$
\end{conj}
Further progress in this direction has been achieved in~\cite{SS:05} and~\cite{GS:17}.
%In our setting of $r=3$ and $j=2$The first open case is to show that $f^{(3)}(n;7,4)=o(n^2)$, which is often referred to as the \defn{$(7,4)$-problem}.

\section{Proof of Theorem~\ref{thm:main}}

\subsection{Upper bound}

It is easy to see that for a given parameter $b\in \bR$, we have
\begin{align}
\max_{\substack{x,y\in \bR\\ \text{s.t. } x\ge 4y, \, x+y\le b}} x+2y= \frac{6}{5}b. \label{optim}
\end{align}
Indeed, assuming $x+y\le b$, we have $5y\le b+4y-x$, and deduce $$x+2y=(x-4y)+6y \le (x-4y)+\frac{6}{5}(b+4y-x) = \frac{6}{5}b -\frac{1}{5}(x-4y)  \le \frac{6}{5}b. $$ Equality holds for $x=4b/5$, $y=b/5$.
Using~\eqref{optim}, we can prove the following lemma.

\begin{lemma}\label{lem:upper}
$f^{(3)}(n;5,3) \le \frac{n^2}{5}$.
\end{lemma}

\proof
Let $G$ be any $\cF^{(3)}(5,3)$-free $3$-graph on $n$~vertices. Clearly, the maximum codegree of $G$ is at most~$2$. For each $i\in\Set{1,2}$, let $G_i$ be the $2$-graph on $V(G)$ whose edges are the pairs $xy$ with $d(xy)=i$. Thus, we have $3e(G)=\sum_{xy}d(xy)=e(G_1)+2e(G_2)$.
The crucial observation is that $e(G_1)\ge 4e(G_2)$. Indeed, for every edge $xy$ in $G_2$, there are distinct $z,z'$ such that $xyz,xyz'\in E(G)$. Note that none of the pairs $xz,yz,xz',yz'$ can be contained in another triple. Thus, $xz,yz,xz',yz'\in E(G_1)$, and none of these pairs is obtained in the same way starting from another edge $x'y'\in E(G_2)$.

Since $e(G_1)+e(G_2)\le \binom{n}{2}$, invoking~\eqref{optim} yields $3e(G)=e(G_1)+2e(G_2)\le \frac{6}{5}\binom{n}{2}$, implying that $e(G)\le \frac{n(n-1)}{5}$, as desired.
\endproof

\subsection{Lower bound}
To establish the lower bound, we use the following well-known result. An \defn{$H$-packing} in a graph $G$ is a collection of edge-disjoint subgraphs of $G$ each isomorphic to~$H$.
\begin{theorem} \label{thm:approx dec}
Let $H$ be any graph and $\eps>0$. For sufficiently large $n$, there exists an $H$-packing in $K_n$ covering all but at most $\eps n^2$ edges of~$K_n$.
\end{theorem}
If an $H$-packing covers all edges of~$G$, it is called an \defn{$H$-decomposition} of~$G$. Wilson~\cite{wilson:76} showed in 1976 that for sufficiently large~$n$, there exists an $H$-decomposition of $K_n$ subject to necessary divisibility conditions. This was recently generalised to hypergraphs in \cite{GKLO:17}. Although one can deduce Theorem~\ref{thm:approx dec} from Wilson's theorem, perhaps the simplest way to prove Theorem~\ref{thm:approx dec} is using a hypergraph matching theorem (cf.~\cite{PS:89}).

We will apply Theorem~\ref{thm:approx dec} to the following special graph.
For $t\in \bN$, define the graph $H_t$ with vertex set $V(H_t)=\Set{a,b,x_1,\dots,x_{2t}}$ and edge set $$E(H_t)=\Set{ab}\cup \set{ax_i,bx_i}{i\in \Set{1,\dots,2t}} \cup \set{x_{2i-1}x_{2i}}{i\in \Set{1,\dots,t}}.$$ Note that
$e(H_t)=5t+1$.
On the same vertex set, we also define the $3$-graph $\hat{H}_t$ with edge set $$E(\hat{H}_t)=\set{ax_{2i-1}x_{2i}, bx_{2i-1}x_{2i}}{i\in \Set{1,\dots,t}}.$$ Hence,
$e(\hat{H}_t)=2t$.
Observe also that every edge of $\hat{H}_t$ is `supported' by a triangle in~$H_t$, that is, whenever $xyz\in E(\hat{H}_t)$, then $xy,xz,yz\in E(H_t)$. In particular, this implies that whenever we are given a collection $\cH$ of edge-disjoint copies of $H_t$ and replace each such copy with a copy of $\hat{H}_t$ on the same vertex set in the obvious way, then the collection of copies of $\hat{H}_t$ is again edge-disjoint, and their union yields a $3$-graph $G$ with $e(\hat{H}_t)\cdot|\cH|$ edges.
Crucially, the $3$-graph $G$ obtained in this way is even $\cF^{(3)}(5,3)$-free. To see this, suppose for a contradiction that $G$ contains three edges $e_1,e_2,e_3$ which span at most five vertices. Clearly, then two of these edges must overlap in two vertices, say $|e_1\cap e_2|=2$. By the above, $e_1$ and $e_2$ cannot arise from different copies of~$H_t$. Consequently, they play the roles of $ax_1x_2$ and $bx_1x_2$, say, in one of the copies $\hat{H}_t'$ of~$\hat{H}_t$.
%Let $H'_t$ denote the underlying copy of $H_t$.
We must also have $|e_3\cap (e_1\cup e_2)|\ge 2$. However, since $H_t[\Set{a,b,x_1,x_2}]$ is complete by construction, $e_3$ must also belong to $\hat{H}_t'$, which yields a contradiction since no such triple exists in~$\hat{H}_t'$.

Observe that for the last step, it is crucial that $ab\in E(H_t)$, as otherwise there might be a triple $e_3$ from another copy of $\hat{H}_t$ which together with $e_1,e_2$ forms a forbidden subgraph. As a result of this construction, the edges which play the role of $ab$ will not be contained in any triple of~$G$. On the other hand, the edges which play the role of one of the edges $x_{2i-1}x_{2i}$ will be contained in two triples of~$G$. By making $t$ large, this can significantly increase the average codegree of~$G$ (and thus the number of edges).

\begin{lemma}\label{lem:lower}
For every $\eps>0$, there exists $n_0\in \bN$ such that for all $n\ge n_0$, there exists an $\cF^{(3)}(5,3)$-free $3$-graph $G_n$ on $n$ vertices with $e(G_n)\ge (\frac{1}{5}-\eps)n^{2}$.
\end{lemma}

\proof
Given $\eps>0$, choose $t\in \bN$ such that $\frac{5t}{5t+1}\ge \frac{1-5\eps}{1-4\eps}$. In the following, we assume that $n$ is sufficiently large. We apply Theorem~\ref{thm:approx dec} to obtain an $H_t$-packing $\cH$ in~$K_n$ such that all but at most $\eps n^2$ edges of $K_n$ are covered. Hence, $e(H_t)|\cH|\ge \binom{n}{2} - \eps n^2$, implying that $|\cH|\ge \frac{1}{5t+1}(\frac{1}{2}-2\eps)n^2$.

Now, define the $3$-graph $G_n$ on $V(K_n)$ as above, by replacing every copy of $H_t$ in $\cH$ with a copy of $\hat{H}_{t}$ in the obvious way. By the above observation, $G_n$ is $\cF^{(3)}(5,3)$-free, and
\begin{align*}
e(G_n) =e(\hat{H}_t) \cdot |\cH| \ge \frac{2t}{5t+1}\left(\frac{1}{2}-2\eps\right)n^2 \ge \left(\frac{1}{5}-\eps\right)n^{2},
\end{align*}
which completes the proof.
\endproof

\bigskip
Clearly, Lemmas~\ref{lem:upper} and~\ref{lem:lower} imply Theorem~\ref{thm:main}.

\section{Further results} \label{sec:further}

As mentioned before, Brown, Erd\H{o}s and S\'os~\cite{BES:73a} suggested that $\lim_{n\to \infty}n^{-2}f^{(3)}(n;6,4)=1/4$. We disprove this by showing the following.
\begin{theorem} \label{thm:64}
$f^{(3)}(n;6,4) \le \frac{3}{14}n^2$.
\end{theorem}

\proof
Let $G$ be any $\cF^{(3)}(6,4)$-free $3$-graph on $n$~vertices. Clearly, $G$ has maximum codegree at most~$3$. It is easy to see that we may assume that $G$ is $\cF^{(3)}(4,3)$-free, as each such subgraph would have to be disconnected from the rest of~$G$.
%
%First, we might argue that we can assume that $G$ is $\cF^{(3)}(k-2,k-3)$-free.
%For this, we claim that any subgraph $F$ of $G$ on $k-2$ vertices with $k-3$ edges is disconnected from the rest of the graph. Indeed, any triple connecting $V(F)$ to $V(G)\sm V(F)$ would add a further edge on the expense of at most $2$ more vertices, yielding a $\cF^{(3)}(k,k-2)$-graph.
%Let $G'$ be the graph obtained by deleting all these components. Let $n'':=n-|V(G')|$. Then $e(G)\le e(G') + n''(k-3) \le \alpha (n-n'')^2 + n''(k-3) \le \alpha n^2$.
%

For each $i\in\Set{1,2,3}$, let $G_i$ be the $2$-graph on $V(G)$ whose edges are the pairs $xy$ with $d(xy)=i$, and let $e_i:=e(G_i)/n^2$. Thus, we have $3e(G)=\sum_{xy}d(xy)=e(G_1)+2e(G_2)+3e(G_3)$ and $e_1+e_2+e_3\le 1/2$.

Let $T_{1}$ be the set of triples $xyz\in E(G)$ with $d(xy)=3$ and $d(xz)=d(yz)=1$. Clearly, we have $|T_{1}|=3e(G_3)$. Moreover, let $T_{2}$ be the set of triples $xyz\in E(G)$ with $d(xy)=d(xz)=2$ and $d(yz)=1$. %Let $t_{2}:=|T_{2}|/n^2$. 
%Clearly, we have $t_{1}=3e_3$. 
Note that $d(xy)+d(xz)+d(yz)=5$ for all $xyz\in T_1\cup T_2$ and $d(xy)+d(xz)+d(yz)\le 4$ for all $xyz\in E(G)\sm (T_1\cup T_2)$.
Double-counting yields
\begin{align*}
e(G_1)+4e(G_2)+9e(G_3) &= \sum_{e\in E(G_1\cup G_2 \cup G_3)}d(e)^2 = \sum_{xyz\in E(G)} d(xy)+d(xz)+d(yz) \\
                          &\le 5 |T_{1}\cup T_{2}| + 4 |E(G)\sm (T_{1}\cup T_{2})| = |T_{1}|+|T_{2}|+4e(G), %= |T_{1}|+|T_{2}| + \frac{4}{3}(e(G_1)+2e(G_2)+3e(G_3))
\end{align*}
which implies $-e(G_1)/3 + 4e(G_2)/3 + 2e(G_3) \le |T_{2}|$.

Moreover, for any pair $xy\in E(G_3)$, by our assumption that $G$ is $\cF^{(3)}(4,3)$-free, there are distinct vertices $z_1,z_2,z_3$ such that $xyz_1,xyz_2,xyz_3\in E(G)$. Let $E_{xy}:=\Set{xz_1,xz_2,xz_3,yz_1,yz_2,yz_3}$. Since $G$ is $\cF^{(3)}(6,4)$-free, we must have $E_{xy}\In E(G_1)$. Similarly, for any triple $xyz\in T_{2}$ with $d(xy)=d(xz)=2$ and $d(yz)=1$, there are distinct vertices $w_1,w_2\in V(G)\sm \Set{x,y,z}$ such that $xyw_1,xzw_2\in E(G)$. Let $E_{xyz}:=\Set{xw_1,yw_1,xw_2,zw_2,yz}$. Clearly, we must have $E_{xyz}\In E(G_1)$. Note that all of the above sets $E_{xy}$ and $E_{xyz}$ are pairwise disjoint. Hence, $e(G_1)\ge 6e(G_3)+5|T_{2}|$. Together with the previous constraint, we obtain $-8e_1/3 + 20e_2/3 + 16e_3 \le 0$.

Using a standard linear programming tool, one can check that
\begin{align*}
\max_{\substack{e_1,e_2,e_3 \ge 0 \\ \text{s.t. } e_1+e_2+e_3\le 1/2  \\ -8e_1/3 + 20e_2/3 + 16e_3 \le 0 } } \frac{1}{3}(e_1+2e_2+3e_3) = \frac{3}{14}.
\end{align*}
This implies $e(G)\le \frac{3}{14}n^2$.
\endproof

We remark that Oleg Pikhurko has improved the constant further to~$\frac{7}{36}$ by considering slightly more complicated configurations.

As mentioned before, the lower bound $1/6$ in~\eqref{best general} is probably not sharp for any $k>4$. 
On the other hand, $1/6$ can be the correct answer when forbidding subgraphs of more than one order larger than~$4$. For instance, we observe the following. It seems plausible that a similar result holds more generally.

\begin{theorem}\label{thm:five and six}
$ex(n;\cF^{(3)}(5,3)\cup \cF^{(3)}(6,4))=(\frac{1}{6}\pm o(1))n^2$.
\end{theorem}

\proof
The lower bound follows from known constructions of $\cF^{(3)}(6,4)$-free Steiner triple systems (and also from the results in~\cite{BW:18,GKLO:18}).
 It remains to show that 
$ex(n;\cF^{(3)}(5,3)\cup \cF^{(3)}(6,4))\le \binom{n}{2}/3$.
Let $G$ be any $(\cF^{(3)}(5,3)\cup \cF^{(3)}(6,4))$-free $3$-graph on $n$ vertices. Clearly, the maximum codegree of $G$ is at most~$2$. For each $i\in\Set{0,1,2}$, let $G_i$ be the $2$-graph on $V(G)$ whose edges are the pairs $xy$ with $d(xy)=i$. Thus, we have $3e(G)=\sum_{xy}d(xy)=e(G_1)+2e(G_2)$.
%The crucial claim is that $e(G_0)\ge e(G_2)$.
We define a map $\phi\colon E(G_2)\to E(G_0)$ as follows. Given an edge $xy\in E(G_2)$, there are unique distinct $z,z'$ such that $xyz,xyz'\in E(G)$. Since $G$ is $\cF^{(3)}(5,3)$-free, we must have $zz'\in E(G_0)$. Let $\phi(xy):=zz'$. Moreover, since $G$ is $\cF^{(3)}(6,4)$-free, $\phi$ must be injective. This implies $e(G_2)\le e(G_0)$ and
thus $3e(G) = e(G_1)+2e(G_2) \le e(G_0)+ e(G_1) +e(G_2)=\binom{n}{2}$.
\endproof

\section{Concluding remarks}

It would be interesting to investigate whether Conjecture~\ref{conj:BES} can be proven without actually determining the limit. For instance, a folklore observation by Katona, Nemetz and Simonovits~\cite{KNS:64} is that for every family of $r$-graphs~$\cF$, a simple averaging argument shows that $\binom{n}{r}^{-1}ex(n;\cF)$ is a decreasing sequence in $[0,1]$, and thus has a limit (called the \defn{Tur\'an density of~$\cF$}). Perhaps similar methods can be used to prove Conjecture~\ref{conj:BES}.
%In fact, one might wonder whether $ex(n;\cF)=\Theta(n^\alpha)$ always implies that $\lim_{n\to \infty}n^{-\alpha} ex(n;\cF)$ exists.

Of course, even if Conjecture~\ref{conj:BES} can be proven in such a way, it would still be desirable to determine the limits. We believe that our methods can be further developed to tackle this. It is probably not too difficult to establish a general upper bound which improves the one in~\eqref{best general}. For instance, a slight adaptation of the proof of Theorem~\ref{thm:64} also yields $f^{(3)}(n;7,5) \le \frac{259}{999}n^2$.  
More work seems needed to improve the lower bound. We believe that our approach of using approximate $H$-decompositions of~$K_n$ to construct $\cF^{(3)}(k,k-2)$-free $3$-graphs with many edges will be useful for general~$k$. However, this requires a stronger decomposition result than Theorem~\ref{thm:approx dec}. For instance, it would be necessary to ensure that in such an approximate $H$-decomposition, forbidden subgraphs are not formed by triples each arising from a different copy of~$H$. (The reason why this extra care was not necessary for $k=5$ is that there are no linear $\cF^{(3)}(5,3)$-graphs.) However, the results in~\cite{BW:18,GKLO:18} give hope that this is possible. 

It would also be interesting to examine the structure of (near-)extremal examples more closely. For some families~$\cF$, there is a unique extremal example $G_n$ with $e(G_n)=ex(n;\cF)$ for all~$n$, and any $\cF$-free $G$ on $n$ vertices with $e(G)\ge (1-o(1))ex(n;\cF)$ must be structurally close to~$G_n$.
For instance, it follows from the proof of Lemma~\ref{lem:upper} that any $\cF^{(3)}(5,3)$-free $3$-graph $G$ on $n$ vertices with $e(G)\ge (1-o(1))n^2/5$ edges has $o(n^2)$ pairs of codegree~$0$, $(2/5\pm o(1))n^2$ pairs of codegree~$1$, and $(1/10\pm o(1))n^2$ pairs of codegree~$2$.

\section*{Acknowledgment}

I am grateful to Felix Joos and Oleg Pikhurko for helpful comments on the manuscript.
%In particular, Oleg Pikhurko's improved bound for the $(6,4)$-problem saved me from proposing a false conjecture.

\bibliographystyle{amsplain_v2.0customized}
\bibliography{References}

\end{document}